\documentclass{amsart}

\usepackage{amsmath,amssymb}

\newtheorem{theorem}{Theorem}
\newcommand{\bt}{\begin{theorem}}
\newcommand{\et}{\end{theorem}} 

\newtheorem{corollary}{Corollary}
\newcommand{\bc}{\begin{corollary}}
\newcommand{\ec}{\end{corollary}} 

\newtheorem{lemma}{Lemma}
\newcommand{\bl}{\begin{lemma}}
\newcommand{\el}{\end{lemma}}

\newtheorem{problem}{Problem}
\newcommand{\bprob}{\begin{problem}}
\newcommand{\eprob}{\end{problem}}

\newcommand{\beq}{\begin{equation}}
\newcommand{\eeq}{\end{equation}}
\newcommand{\F}{\ensuremath{ \mathbf F }}
\newcommand{\N}{\ensuremath{ \mathbf N }}
\newcommand{\Z}{\ensuremath{\mathbf Z}}

\newcommand{\R}{\ensuremath{\mathbf R}}
\newcommand{\benum}{\begin{enumerate}}
\newcommand{\eenum}{\end{enumerate}}

\DeclareMathOperator{\qqand}{\qquad\text{and}\qquad}

\title{Sidon sets with $\Delta$-separated sumsets}

\author{Melvyn B. Nathanson}

\address{Department of Mathematics, Lehman College (CUNY), Bronx, NY, USA}  
\email{melvyn.nathanson@lehman.cuny.edu}

\date{\today}


\begin{document}

\begin{abstract}
The set $A$ is a $B_h$-set 
if every element of the sumset $hA$ has a unique representation as a sum of $h$ 
elements of $A$.  A $B_2$-set is also called a Sidon set.  
A $B_{h,\Delta}$-set is a $B_h$-set $A$ whose sumset $hA$ is $\Delta$-separated, 
that is, $x'-x\geq \Delta$ for all $x,x' \in hA$ with $x <x'$.  
Upper and lower bounds are obtained for the cardinality of the 
largest $B_{2,\Delta}$-sets contained in   $\{1,2,\ldots, n\}$, 
that is, sets $A \subseteq \{1,2,\ldots, n\}$ such that, if $a,b,c,d \in A$ 
and $\{a,b\} \neq \{c,d\}$, then $|(a+b)-(c+d)| \geq \Delta$.  
\end{abstract}

\maketitle 

\section{Sidon sets}

 A  set $X$ of real numbers  is \emph{$\Delta$-separated} if $\Delta > 0$ and 
\[
x'-x\geq \Delta
\]
for all $x,x' \in X$ with $x < x'$. 
Every set of integers is 1-separated. 
Every translate of a $\Delta$-separated set is $\Delta$-separated. 

Let $A$ be a nonempty subset of an additive abelian group or semigroup $X$.   
The \emph{representation function} 
$r_{A,h}(x)$ counts the number of $h$-tuples $(a_1,\ldots, a_h)$ 
of elements of $A$ such that 
\[
a_1 \leq \cdots \leq a_h 
\]
and
\[
a_1+\cdots + a_h = x.
\]
The set $A$  is called a \emph{$B_h[g]$-set} 
if every  element $x \in X$ has at most $g$ representations as the sum of $h$ elements of $A$, 
that is, if $r_{A,h}(m) \leq g$ for all $m \in \Z$.  
A $B_h[1]$-set  is called a \emph{$B_h$-set} and a $B_2$-set is called a \emph{Sidon set}. 
A Sidon set satisfies the following condition:  
For all $a,b,c,d \in A$, we have $a+b = c+d$ if and only if 
$\{a,b\} = \{c,d\}$. 
Two classical problems in additive number theory are to construct large  
infinite $B_h[g]$-sets of positive integers and to compute the size of the largest 
$B_h[g]$-set A contained in the integer interval $\{1,2,\ldots, n\}$.

Let $A$ be a $B_h$-set  of integers.  For all $a_1,\ldots, a_h, a'_1,\ldots, a'_h \in A$ with 
\[
a_1 \leq \cdots \leq a_h \qqand a'_1 \leq \cdots \leq a'_h
\]
we have 
\[
a_1 + \cdots + a_h = a'_1 + \cdots + a'_h
\]
if and only if 
\[
(a_1, \ldots, a_h) = (a'_1, \ldots, a'_h). 
\]
Equivalently, if $A$ is a $B_h$-set of integers, then conditions 
\[
a_1 \leq \cdots \leq a_h \qqand a'_1 \leq \cdots \leq a'_h
\]
and 
\[
(a_1, \ldots, a_h) \neq (a'_1, \ldots, a'_h) 
\]
imply 
\[
|(a_1 + \cdots + a_h) - (a'_1 + \cdots + a'_h)| \geq 1.
\]
For every positive integer $\Delta$, it is natural to consider $B_h$-sets 
whose $h$-fold sumsets are $\Delta$-separated, 
that is, sets with the property that the conditions 
\[
a_1 \leq \cdots \leq a_h \qqand a'_1 \leq \cdots \leq a'_h
\]
and 
\[
(a_1, \ldots, a_h) \neq (a'_1, \ldots, a'_h) 
\]
imply  
\[
|(a_1 + \cdots + a_h) - (a'_1 + \cdots + a'_h)| \geq \Delta. 
\]
We call  such sets $B_{h,\Delta}$-sets.  
A $B_{h,\Delta}[g]$-set is a $B_h[g]$-set whose sumset $hA$ is $\Delta$-separated. 
Every translate of a $B_{h,\Delta}[g]$-set is a $B_{h,\Delta}[g]$-set.


\section{$\Delta$-separated sets}

The \emph{$\Delta$-dilation} of the set $X$ is  
\[
\Delta \ast X = \{\Delta x: x \in X\}. 
\]
If $\Delta' > 0$ and $X$ is $\Delta'$-separated, 
then the set $\Delta\ast X$ is $\Delta\Delta'$-separated. 
Because every set $X$ of integers is 1-separated, 
for every positive integer $\Delta$, the set $\Delta\ast X$ is a $\Delta$-separated 
set of integers.

\bl                    \label{Delta:lemma:Delta-size} 
Let $A$ be a  $\Delta$-separated set of real numbers.  
If $A \subseteq [u,v]$, then 
\[
|A| \leq 1 + \frac{v-u}{\Delta}.
\]

\el

\begin{proof}
Let $A = \{a_1,\ldots, a_k\}$, where $|A|=k$ and 
\[
u \leq a_1 < a_2 < \cdots < a_k \leq v.
\]
Then 
\[
v-u \geq a_k-a_1 = \sum_{i=1}^{k-1} (a_{i+1}-a_i ) \geq (k-1)\Delta
\]
and so 
\[
k \leq 1 + \frac{v-u}{\Delta}.
\]
This completes the proof. 
\end{proof}

Let $A$ be a set of integers. 
The \emph{counting function} $A(x)$ counts the number of positive integers in $A$ 
that do not exceed $x$: 
\[
A(x) = \sum_{\substack{a \in A \\ 0 < a \leq x}} 1.
\]
We have $0 < a \leq x$ if and only if $0 < \Delta a \leq \Delta x$.   
Moreover, if $a' \in \Delta\ast A$, then $a' \equiv 0 \pmod{\Delta}$ 
and so $A(x) = (\Delta \ast A)(\Delta x)$.  Equivalently, 
\beq          \label{Delta:count}
(\Delta \ast A)(x) = A\left( \frac{x}{\Delta} \right) 
\eeq
for all $x > 0$. 

\bl           \label{Delta:lemma:reduce}
Let $A$ be a nonempty set of real numbers and let $h \geq 2$ and $\Delta > 0$.  
If the sumset $hA$ is $\Delta$-separated, then the sumset $(h-1)A$ 
is $\Delta$-separated. 
The converse is not necessarily true.    
\el

\begin{proof} 
If the sumset $(h-1)A$ is not $\Delta$-separated, then there exist 
$a,b \in (h-1)A$ with $a \neq b$ 
and $|a-b| < \Delta$.  For any $c \in A$, the integers $a+c \in hA$ and 
$b+c \in hA$ satisfy $a+c \neq b+c$  and 
\[
|(a+c)-(b+c)| = | a-b | < \Delta
\]
and so the sumset $hA$ is not $\Delta$-separated.  
Thus, if $hA$ is $\Delta$-separated, then $(h-1)A$ is $\Delta$-separated. 

To prove that the converse is not necessarily true, 
let $\Delta  > 1$ and  choose integers $n, u,v$ 
with $n > \Delta$ and $v-u> 2n + \Delta$. 
The set $A = \{u, u+n,  v - n, v+1\}$ is $\Delta$-separated but the sumset $2A$ 
is not $\Delta$-separated because $a = u+ (v+1 ) \in 2A$, $b = u+v = (u+n) + (v-n) \in 2A$, 
and $|a-b| = 1 < \Delta$. 
This completes the proof. 
\end{proof}

\bl                            \label{Delta:lemma:dilateBh}
Let $A$ be a nonempty set of positive integers.  
For all positive integers  $\Delta$ and $g$, 
 the set $A$ is a $B_h[g]$-set
if and only if the dilation $\Delta\ast A$ is a $B_{h,\Delta}[g]$-set. 
\el

\begin{proof} 
Let $A$ be a set of integers.  
The dilated set $\Delta \ast A$ is $\Delta$-separated.  
For all $a_1,\ldots, a_h, b_1,\ldots, b_h \in A$, we have 
\[
\sum_{i=1}^h \Delta a_i = \sum_{i=1}^h \Delta b_i 
\]
if and only if 
\[
\Delta \sum_{i=1}^h a_i = \Delta  \sum_{i=1}^h b_i 
\]
if and only if 
\[
\sum_{i=1}^h a_i = \sum_{i=1}^h b_i 
\]
and so $r_{A,h}(m) = r_{\Delta\ast A,h}(\Delta m)$ for all integers $m$. 
Moreover, $r_{\Delta\ast A,h}(m') = 0$ if $m' \not\equiv 0 \pmod{\Delta}$, 
and so $r_{A,h}(m) \leq g$ for all integers $m$ if and only if 
$ r_{\Delta\ast A,h}(m) \leq g$  for all integers $m$. 
This completes the proof. 
\end{proof}

\section{ A lower bound for $B_h$-sets with $\Delta$-separated sumsets}

Let $F_{h,\Delta}(n)$ denote the size of the largest $B_{h,\Delta}$-set contained in
$\{1,2,\ldots, n\}$.

 \bt                        \label{Delta:theorem:LowerBound}
 There exists $0 < \theta< 1$ such that, 
for all  positive integers $h$ and $\Delta$ and for all sufficiently large  $n$, 
\beq                    \label{Delta:LowerBound}
 F_{h,\Delta}(n)   \geq  \left(\frac{n}{\Delta}\right)^{1/h} 
 -   \left(\frac{n}{\Delta}\right)^{\theta/h}  
\eeq
 and 
\beq                    \label{Delta:LowerBound-inf}
 \liminf_{n\rightarrow \infty} \frac{  F_{h,\Delta}(n)  }{  \left(n/\Delta\right)^{1/h} } \geq 1.
 \eeq
 \et
 
 \begin{proof}
 We start with a  classical argument of Bose-Chowla~\cite{bose-chow62}, 
 which extended earlier results of Bose~\cite{bose42}, Chowla~\cite{chow44}, 
and Singer~\cite{sing38}. 

Let $q$ be a prime power and consider the finite fields $\F_{q}$ and $\F_{q^h}$.  
Let $\alpha$ be a primitive element for the finite field extension $\F_{q^h}$ 
over $\F_q$.  Then $\alpha \notin \F_q$ 
and $\F_{q^h} \setminus \{0\} = \{\alpha^i: i =1, \ldots, q^h -1 \}$. 

Let $\F_q   = \{\lambda_1,\ldots, \lambda_q\}$. 
Because $\alpha \notin \F_q$, for all integers $i \in [1,q]$ 
we have $\alpha - \lambda_i \neq 0$ 
  and so  $\alpha - \lambda_i \in \F_{q^h}  \setminus \{0\} $.   
 Let $a_i$ be the unique integer in $[1, q^h -1]$ such that 
 \[
\alpha - \lambda_i =  \alpha^{a_i}.
 \]
Then 
\[
A = \{a_i:i  = 1,\ldots, q \} \subseteq [1, q^h -1]
\]
and $|A| = q$.  We shall prove that $A$ is a $B_h$-set. 

Let $(a_{u_1},\ldots, a_{u_h})$ and 
$(a_{v_1},\ldots, a_{v_h})$ be increasing sequences in $ A$ 
such that 
\[
 a_{u_1} + \cdots + a_{u_h} =  a_{v_1}  + \cdots +  a_{v_h}.
\]
Then 
\begin{align*}
\prod_{i=1}^h  (\alpha - \lambda_{u_i}) 
& = \prod_{i=1}^h  \alpha^{a_{u_i}} =  \alpha^{\sum_{i=1}^h a_{u_i}} \\
& =  \alpha^{\sum_{i=1}^h a_{v_i}} = \prod_{i=1}^h  \alpha^{a_{v_i}}\\
& = \prod_{i=1}^h  (\alpha - \lambda_{v_i}).
\end{align*} 
The polynomials 
\[
f(x)  = \prod_{i=1}^h  (x - \lambda_{u_i}) - x^h  
\]
and 
\[
g(x) = \prod_{i=1}^h  (x - \lambda_{v_i}) - x^h 
\]
are nonzero and have degrees at most $h-1$ and 
\begin{align*}
f(\alpha) & = \prod_{i=1}^h  (\alpha - \lambda_{u_i}) - \alpha^h  = 
\prod_{i=1}^h  (\alpha - \lambda_{v_i}) - \alpha^h = g(\alpha). 
\end{align*} 
Because the minimal polynomial of $\alpha$ in $\F_q[x]$ has degree $h$ 
and because $ f(x) - g(x)$ has degree at most $h-1$ 
and $f(\alpha) - g(\alpha) = 0$, it follows that $f(x) = g(x)$.   
Equivalently, 
\[
 \prod_{i=1}^h  (x- \lambda_{u_i}) = \prod_{i=1}^h  (x- \lambda_{v_i}).
\]
This implies $\{\lambda_{u_1},\ldots, \lambda_{u_h} \} 
= \{ \lambda_{v_1},\ldots, \lambda_{v_h} \}$ as multisets 
and so $(a_{u_1},\ldots, a_{u_h}) = (a_{v_1},\ldots, a_{v_h})$. 
Thus, for every prime power $q$, there is a   $B_h$-set of size $q$ contained in $\{1,2,\ldots, q^h -1\}$.  
In particular, for every prime $p$, there is a $B_h$-set $A$ of size $p$ 
contained in $\{1,2,\ldots, p^h -1\}$. 
By Lemma~\ref{Delta:lemma:dilateBh}, the set $A' = \Delta\ast A$ 
is a $B_{h,\Delta}$-set contained in 
$\{1, 2, \ldots, \Delta  p^h  \}$ with 
\[
|A'| = |A| = p. 
\] 

From prime number theory, we know that there exists $0 < \theta < 1$ 
such that, for all sufficiently large integers $n$ and 
\[
m = \left(\frac{n}{\Delta}\right)^{1/h}
\]
there is a prime $p$  with 
\beq \label{Delta:li}
m-m^{\theta} \leq p \leq m.  
\eeq
(Runbo Li~\cite{li26}  proved~\eqref{Delta:li} with $\theta = 0.52$.) 
Because 
\[
 \Delta  p^h \leq  \Delta  m^h = n 
\]
and 
\[
A' \subseteq \{1,\ldots, \Delta p^h \} \subseteq \{1,\ldots, n\} 
\]
we have  
\[ 
F_{h,\Delta}(n)  \geq |A'|  = p  \geq m-m^{\theta} = \left(\frac{n}{\Delta}\right)^{1/h} - \left(\frac{n}{\Delta}\right)^{\theta /h} . 
\] 
This immediately implies 
\[
 \liminf_{n\rightarrow \infty} \frac{  F_{h,\Delta}(n)  }{  \left(n/\Delta\right)^{1/h} } \geq 1 
\]
and completes the proof.  
\end{proof}

\section{ An upper bound for $B_2$-sets with $\Delta$-separated sumsets} 
There is an elementary upper bound for $ F_{h,\Delta}[g](n)$. 

\bt             \label{Delta:theorem:UpperBound-easy}
For all positive integers $g$ and $h$, 
 \beq             \label{Delta:UpperBound-easy}
 \limsup_{n\rightarrow \infty} \frac{  F_{h,\Delta}[g](n)  }{  \left(gn/\Delta\right)^{1/h} }
 \leq (h!h)^{1/h}  < \infty.
\eeq
\et

\begin{proof}
Let $A$ be a $B_{h,\Delta}[g]$-set contained in $\{1,2,\ldots, n\}$ 
with $|A| = k$.  Then $hA \subseteq [h,hn]$.  
Because the sumset $hA$ is $\Delta$-separated, 
Lemma~\ref{Delta:lemma:Delta-size} implies that 
\[
|hA| \leq 1 + \frac{hn-h}{\Delta}. 
\]
Let $r_{A,h}(m)$ be the $h$-fold representation function of $A$.  
Because $A$ is a $B_{h,\Delta}[g]$-set, we have $r_{A,h}(m) \leq g$ 
for all integers $m$.   Because $|A|= k$, the number of $h$-tuples 
$(a_1,\ldots, a_h) \in A^h$ with $a_1\leq \cdots \leq a_h$ is 
$\binom{h+k-1}{h}$ and so
\begin{align*}
\frac{k^h}{h!} 
&  \leq \binom{h+k-1}{h} = \sum_{m=h}^{hn} r_{A,h}(m) = \sum_{m\in hA} r_{A,h}(m) \\
& \leq g|hA| \leq g\left( 1 + \frac{hn-h}{\Delta}\right) \\
& < \frac{ghn}{\Delta} + g.   
\end{align*}
We  obtain  
\[
k < \left(   \frac{h!ghn}{\Delta} +  h!g  \right)^{1/h}.
\]
Inequality~\eqref{Delta:UpperBound-easy} follows immediately. 
This completes the proof. 
\end{proof}

The \emph{positive difference set} of the set $A$ is 
\[
D^+(A) = \{a-b: a,b \in A  \text{ and } a > b\}.
\]

\bl               \label{Delta:lemma:reduce-diff}
If  the set $A$ has a $\Delta$-separated positive difference set, 
then $A$ is $\Delta$-separated.
\el

\begin{proof}
Suppose that  $A$ has a $\Delta$-separated  positive difference set.  
If $|A|=1$ or 2,  then $A$ is $\Delta$-separated.  
Let $|A| \geq 3$. 
If $A$ is not $\Delta$-separated, then there exist $a,b \in A$ with $0 < a-b < \Delta$. 
Choose $c \in A\setminus \{ a,b\}$. 
If $c < b <a$, then 
\[
0 < (a-c)-(b-c) = a-b < \Delta 
\]
 which is absurd. 
If $b < a < c$, then 
\[
0 < (c-b)-(c-a) = a-b < \Delta 
\]
which is absurd.  
 If $b< c < a$, then 
\[
0 < (a-b)-(c-b) = a-c < a-b < \Delta
\]
which is absurd. 
This completes the proof. 
\end{proof}

The set $A$ has a \emph{unique positive difference set} 
if, for all $a,b,c,d \in A$ with $a>b$ and $c > d$, 
we have $a-b=c-d$ if and only if $a=c$ and $b=d$.  
The set $A$ has a \emph{$\Delta$-separated unique  positive difference set} 
if, for all $a,b,c,d \in A$ with $a>b$ and $c > d$, 
we have 
\[
|(a-b)-(c-d)| < \Delta 
\]
if and only if $a=c$ and $b=d$.

There is a simple relation between  
$B_{2,\Delta}$-sets and sets with a $\Delta$-separated unique positive difference set.

\bl                     \label{Delta:lemma:Delta-difference}
The set $A$ is a $B_{2,\Delta}$-set if and only if $A$  is a set 
with a $\Delta$-separated unique positive  difference set.
\el

\begin{proof} 
Let $A$ be a $B_{2,\Delta}$-set. 
 By Lemma~\ref{Delta:lemma:reduce}, the set $A$ is also $\Delta$-separated.  
Let $a,b,c,d \in A$ satisfy $a > b$ and $c > d$ and 
\[
|(a-b)-(c-d)|< \Delta.   
\]
Then 
\[
|(a+d)-(b+c)| < \Delta.   
\]
Because $A$ is a $B_{2,\Delta}$-set, 
we have $\{a,d\}=\{b,c\}$. Then $a> b$ implies $a=c$ and $b=d$ and so 
$A$ has a  $\Delta$-separated unique difference set.  

Conversely, let $A$ have a $\Delta$-separated unique positive  difference set. 
 By Lemma~\ref{Delta:lemma:reduce-diff}, the set $A$ is also $\Delta$-separated.  
If $a,b,c,d \in A$ satisfy 
\[
|(a+b)-(c+d)|< \Delta  
\]
then  
\[
|(a-d)-(c-b)| < \Delta. 
\]
If $a=d$, then $|c-b| < \Delta$ and so, because $A$ is $\Delta$-separated, 
$b=c$ and $\{a,b\}=\{c,d\}$.  

If $a > d$, then $a-d \geq \Delta$  because $A$  is a $\Delta$-separated. 
If $c \leq b$, then 
\[
\Delta \leq  a-d \leq (a-d) + (b-c) = |(a-d)-(c-b)| < \Delta
\]
which is absurd.  
Therefore, $c > b$. 
Because $A$ has a $\Delta$-separated unique difference set, we have 
$a=c$ and $b=d$ and so $\{a,b\} = \{c,d\}$. 
Thus, $A$ is a $B_{2,\Delta}$-set. 

If $a < d$, then apply the argument to $|(d-a)-(b-c)| < \Delta$. 
This completes the proof. 
\end{proof}

Erd\H os and Tur\' an~\cite{erdo-tura41} proved that, for every positive integer $n$,  
if $A$ is a $B_2$-set contained in $\{1,2,\ldots, n\}$, then 
\[
|A| < n^{1/2} + O\left(n^{1/4}\right). 
\] 
Equivalently, 
\[
F_2(n) < n^{1/2} + O\left(n^{1/4}\right). 
\] 
Lindstr\" om~\cite{lind69} proved that 
\[
F_2(n) <n^{1/2} + n^{1/4} + 1.
\] 
There is a recent improvement by 
Balogh, F\"uredi, and Roy~\cite{balo23}: 
\[
F_2(n) <n^{1/2} + 0.998 n^{1/4} + o(1).
\] 
 
In an early version of this paper~\cite{nath26}, I applied the 
 Erd\H os-Tur\' an argument to prove that 
\[
F_{2,\Delta}(n) < \left(\frac{2n}{\Delta}\right)^{1/2} + O\left(n^{1/4}\right). 
\] 
Kevin O'Bryant (personal communication) observed 
that one could use Lindstr\" om's method to obtain the following better result. 

\bt         \label{Delta:theorem:UpperBound}
For all positive integers $\Delta$, 
\beq                    \label{Delta:UpperBound}
F_{2,\Delta}(n) \leq  \left(\frac{n}{\Delta}\right)^{1/2} + O\left(n^{1/4}\right) 
\eeq
and
\beq                    \label{Delta:UpperBound-sup}
 \limsup_{n\rightarrow \infty} \frac{  F_{2,\Delta}(n)  }{  \left(n/\Delta\right)^{1/2} } \leq 1.
 \eeq
\et

\begin{proof} 
Let $\Delta$ and $n$ be positive integers with $n \geq \Delta+2$ 
and let $A$ be a $B_{2,\Delta}$-set contained in $\{1,2,\ldots, n\}$ with $|A| = k$. 
Following the classical argument of Lindstr\" om~\cite{lind69}, 
we let  $A = \{a_1,a_2,\ldots, a_k\}$ with $a_1 < a_2 < \cdots < a_k$ 
and consider the positive difference set 
\[
D^+(A) = \{a_v-a_u: 1 \leq u < v \leq k\}. 
\]
For all $j \in [1,k-1]$, a ``level $j$ difference'' 
is a positive integer of the form  $a_{i+j}-a_i$ for some $i \in [1,k-j]$.  
By Lemma~\ref{Delta:lemma:Delta-difference}, 
the $B_{2,\Delta}$-set $A$ has a $\Delta$-separated unique positive 
difference set $D^+(A)$, and so $D^+(A)$  contains exactly 
$k-j$ distinct level $j$ differences. 
The sum of these $k-j$ ``level $j$ differences'' is 
\begin{align*} 
S_j & = \sum_{i=1}^{k-j} (a_{i+j} - a_i) =  \sum_{i=j+1}^k a_i - \sum_{i=1}^{k-j}  a_i. 
\end{align*} 
If $k-j \leq j$, then 
\[
S_j  = \sum_{i=j+1}^k a_i - \sum_{i=1}^{k-j}  a_i  
\leq  \sum_{i=j+1}^k a_i  \leq (k-j)a_k \leq jn. 
\]
If $k-j \geq j +1$, then 
\[
S_j  = \sum_{i=j+1}^k a_i - \sum_{i=1}^{k-j}  a_i 
= \sum_{i= k-j+1}^k a_i - \sum_{i=1}^j a_i  \leq ja_k \leq jn. 
\]
For all ${\ell} \in [1,k-1]$,  
\beq                                \label{Delta:Tm-upper}
T_{\ell} = \sum_{j=1}^{{\ell}} S_j \leq \sum_{j=1}^{{\ell}} jn =  \frac{{\ell}({\ell}+1)n}{2}  
\eeq
is the sum of the positive integers in the set $D^+(A)$ that are differences 
of level at most ${\ell}$.  

Because there are exactly $k-j$  differences of level $j$ 
for all $j \in [1,k-1]$, there are exactly 
\[
\sum_{j=1}^{\ell} (k-j)= k{\ell} - \frac{{\ell}({\ell}+1)}{2} = {\ell}\left( k-\frac{{\ell}+1}{2} \right) = {\ell}m
\]
positive integers in $D^+(A)$ that are differences 
of level at most ${\ell}$, where 
\[
m = k-\frac{{\ell}+1}{2}. 
\]
These integers form a $\Delta$-separated set of size ${\ell}m$.  
 The $\Delta$-separated set of size ${\ell}m$ with the smallest positive integers 
 is the arithmetic progression 
 \[
 \{1+i\Delta: i = 0,1,\ldots, {\ell}m-1\}
 \]
  and so 
\beq                                \label{Delta:Tm-lower}
T_{\ell} \geq  \sum_{i=0}^{{\ell}m-1} (1+i\Delta) = {\ell}m + \frac{({\ell}m-1){\ell}m \Delta}{2} 
>  \frac{({\ell}m-1){\ell}m \Delta}{2}.
\eeq
Combining inequalities~\eqref{Delta:Tm-upper} and~\eqref{Delta:Tm-lower}, 
we obtain 
\[
   \frac{({\ell}m-1){\ell}m \Delta}{2} < \frac{{\ell}({\ell}+1)n}{2} 
 \]
and so 
\[
{\ell}\Delta m^2 -\Delta m = ({\ell}m-1)m \Delta < ({\ell}+1)n. 
\]
Dividing by ${\ell}\Delta$ and completing the square gives 
\[
\left( m - \frac{1}{2{\ell}}\right)^2 = m^2  -\frac{m}{{\ell}}  + \frac{1}{4{\ell}^2} < \left(1+\frac{1}{{\ell}}\right) \frac{n}{\Delta} + \frac{1}{4{\ell}^2}. 
\]
From the inequalties $\sqrt{x+y} < \sqrt{x}  + \sqrt{y}$ and $\sqrt{1+x}  < 1 + x/2$ 
for $x > 0$ and $y > 0$, we obtain   
\begin{align*} 
m - \frac{1}{2{\ell}} 
& < \sqrt{ \left(1+\frac{1}{{\ell}}\right) \frac{n}{\Delta} + \frac{1}{4{\ell}^2} } \\ 
& \leq \left(1+\frac{1}{2{\ell}}\right) \sqrt{\frac{n}{\Delta} } + \frac{1}{2{\ell}} \\ 
\end{align*}
and so 
\[
|A| = k = m + \frac{{\ell}+1}{2} < \sqrt{\frac{n}{\Delta} } + \frac{1}{2{\ell}}\sqrt{\frac{n}{\Delta} }  + \frac{{\ell}+1}{2} + \frac{1}{{\ell}}.
\]
Choosing ${\ell} = \left[n^{1/4}\right]$ gives 
\[
|A|  < \left(\frac{n}{\Delta}\right)^{1/2} + O\left(n^{1/4}\right). 
\]
This immediately implies 
\[
 \limsup_{n\rightarrow \infty} \frac{  F_{2,\Delta}(n)  }{  \left(n/\Delta\right)^{1/2} } \leq 1 
\]
and completes the proof.  
\end{proof} 

\bt 
\[
F_{2,\Delta}(n) \sim \left(\frac{n}{\Delta}\right)^{1/2}.
\]
\et

\begin{proof} 
This follows from Theorem~\ref{Delta:theorem:LowerBound} with $h=2$ 
and from Theorem~\ref{Delta:theorem:UpperBound}.
\end{proof}

\section{Open problems}
\subsection{$B_h$-sets of integers} 
\bprob
A general problem is to extend known results about $B_h$-sets to $B_{h,\Delta}$-sets 
  and known results about $B_h[g]$-sets to $B_{h,\Delta}[g]$-sets. 
For example, can we improve Theorems~\ref{Delta:theorem:UpperBound-easy} 
and~\ref{Delta:theorem:UpperBound}?  
For integers $h \geq 3$, determine if the limits 
\[
\lim_{n\rightarrow \infty} \frac{F_{h,\Delta}(n)}{(n/\Delta)^{1/h} }
\]
exist.    This is unknown even in the classical case $\Delta =1$.   
Compute the limits if they do exist.  

Let $g\geq 2$ and  $h \geq 2$.  Determine if the limits 
\[
\lim_{n\rightarrow \infty} \frac{F_{h,\Delta}[g](n)}{(n/\Delta)^{1/h} }
\]
exist and compute them if they exist.  
\eprob

\bprob
For $h \geq 2$, how many $B_{h,\Delta}$-sets of maximum size are contained in 
$\{1,2,\ldots, n\}$? 
Classify the different possible structures 
of $B_{h,\Delta}$-sets of maximum size that are contained in $\{1,2,\ldots, n\}$. 
\eprob

\subsection{$B_h$-sets in semigroups} 
This is a new class of problems. 
A $B_h$-set in an additive abelian group or semigroup $G$ is a subset $A$ of $G$ such that,
for all $h$-tuples $(a_1,\ldots, a_h)$ and $(a'_1,\ldots, a'_h)$ of elements of $A$, 
we have 
\[
a_1 + \cdots + a_h = a'_1 + \cdots + a'_h
\]
if and only if there is a permutation $\sigma$ of $\{1,2,\ldots, h\}$ such that 
$a'_i = a_{\sigma(i)}$ for all $i \in \{1,2,\ldots, h\}$. 
There is the analogous definition of $B_h[g]$-sets for all $g \geq 1$.

\bprob
For positive integers $h$ and $g$, 
compute the maximum size of a $B_h[g]$-set in the finite field $\F_p = \Z/p\Z$. 
\eprob

\bprob
In the group $\Z^d$,  let $m_1,\ldots, m_d$ be positive integers 
and let $P(m_1,\ldots, m_d)$ 
be the parallelpiped of lattice points $(a_1,\ldots, a_d)$ with 
$1 \leq a_j\leq m_j$ for all $j \in \{1,\ldots, d\}$. 
For positive integers $h$ and $g$, 
compute the maximum size of a $B_h[g]$-set in the parallelepiped $P(m_1,\ldots, m_d)$ 
or in $\Gamma \cap \Z^n$, where $\Gamma$ is a convex subset of $\R^n$.  
\eprob

Let $G$ be an additive abelian group on which there is an ``absolute value,'' 
denoted $| \cdot |$, such that 
\[
|0|=0 
\]
and 
\[
|x|  = |-x| > 0 \text{ for all } x \in G\setminus \{0\}.
\]
(We do not assume a triangle inequality.) 
Let $\Delta > 0$.  A subset $X$ of $G$ 
is $\Delta$-separated if $|x-x'| \geq \Delta$ for all $x,x' \in A$ with $x \neq x'$. 
A $B_{h,\Delta}$-set is a $B_h$-set  $A$ in $G$ whose $h$-fold sumset  
$hA$ is $\Delta$-separated. 

There are many absolute values on a group.  
If we define $|0|=0$ and $|x| = \Delta$ for all $x \in G\setminus \{0\}$, then 
every subset of $G$ is $\Delta$-separated and every $B_h[g]$-set is a 
$B_{h,\Delta}[g]$-set.  If we define $|0|=0$ and $|x| = \Delta/2$ 
for all $x \in G\setminus \{0\}$, then 
no subset of $G$ is $\Delta$-separated and there is no $B_{h,\Delta}[g]$-set in $G$.  
We consider $B_{h,\Delta}[g]$-sets with respect to ``interesting'' absolute values.  
Here are two examples. 

\bprob
Consider the group $\Z^d$ with the usual absolute value.  
Determine the maximum size of a $B_{h,\Delta}$-set contained in the 
parallelepiped $P(m_1,\ldots, m_d)$ 
or in $\Gamma \cap \Z^n$, where $\Gamma$ is a convex subset of $\R^n$.  
\eprob

\bprob
Let $p$ be an odd prime.  
Define an absolute value in the finite field $\Z/p\Z$ as follows:
For all congruence classes $a + p\Z$, let $|a+p\Z| = c$, where $c$ is the 
unique integer in $\{0,1,2,\ldots, (p-1)/2\}$ such that $a \equiv \pm c \pmod{p}$. 
Determine the maximum size of a $B_{h,\Delta}$-set contained in $\Z/p\Z$. 
\eprob 

\subsection{$B_h$-sets in sets with bounded gaps} 
Consider the pair $(r,s)$ with $0 \leq s < r$.
The set $X = \{x_i\}_{i=1}^n$ be a set of real numbers 
with $x_i < x_{i+1}$ for all $i = 1,\ldots, n-1$.  
The set $X$  is \emph{$(r,s)$-bounded} if 
\[
s < x_{i+1}-x_i \leq r 
\]
for all $i = 1,\ldots, n-1$.  
Every translate of an $(r,s)$-bounded set is $(r,s)$-bounded.  
Consideration of  $(r,s)$-bounded  sets has led to new problems 
and results about van der Waerden numbers (Fox-Schildkraut~\cite{fox-schi-26}) 
and, similarly, suggests new problems about $B_h$-sets.  

\bprob
What is the size of the largest $B_h$-set 
(or $B_{h,\Delta}$-set or $B_{h,\Delta}[g]$-set) contained in an $(r,s)$-bounded 
set of integers of size $n$? 
\eprob

\bprob
Let $X = \{x_1,\ldots, x_n\}$ be a set of  integers such that 
$x_i \in [(i-1)r+1,ir]$ for all $i = 1,\ldots, n$.  (Such sets were consider in 
Nathanson~\cite{nath80}.)  
 What is the size of the largest $B_h$-set 
(or $B_{h,\Delta}$-set or $B_{h,\Delta}[g]$-set) contained in $X$? 
\eprob

\def\cprime{$'$} \def\cprime{$'$}
\providecommand{\bysame}{\leavevmode\hbox to3em{\hrulefill}\thinspace}
\providecommand{\MR}{\relax\ifhmode\unskip\space\fi MR }
\providecommand{\MRhref}[2]{
  \href{http://www.ams.org/mathscinet-getitem?mr=#1}{#2}
}
\providecommand{\href}[2]{#2}

\end{document}